\documentclass[a4, 10pt]{amsart}
\usepackage{amssymb}
\usepackage{amstext}
\usepackage{amsmath}
\usepackage{amscd}
\usepackage{latexsym}
\usepackage{amsfonts}

\theoremstyle{plain}
\newtheorem{thm}{Theorem}[section]
\newtheorem*{thm*}{Theorem}
\newtheorem*{cor*}{Corollary}

\newtheorem{lem}[thm]{Lemma}

\newtheorem{claim}{Claim}
\newtheorem*{claim*}{Claim}

\theoremstyle{definition}

\theoremstyle{remark}
\newtheorem*{pf}{{\sl Proof}}

\newtheorem*{tpf}{{\sl Proof of Theorem \ref{main}}}

\newtheorem*{cpf}{{\sl Proof of Claim}}

\numberwithin{equation}{thm}

\def\Ext{\mathrm{Ext}}

\def\p{\mathfrak p}
\def\q{\mathfrak q}

\def\P{\mathfrak P}

\def\depth{\mathrm{depth}}
\def\Supp{\mathrm{Supp}}
\def\Ann{\mathrm{Ann}}

\def\id{\mathrm{id}}

\def\height{\mathrm{ht}}

\def\Spec{\mathrm{Spec}}

\def\xx{\text{\boldmath $x$}}

\def\CM{\mathrm{CM}}

\begin{document}

\title[Modules locally of finite injective dimension]{A characterization of modules locally of finite injective dimension}
\author{Ryo Takahashi}
\address{Department of Mathematics, School of Science and Technology, Meiji University, 1-1-1 Higashimita, Tama-ku, Kawasaki 214-8571, Japan}
\email{takahasi@math.meiji.ac.jp}
\begin{abstract}
In this note, we characterize finite modules locally of finite injective dimension over commutative Noetherian rings in terms of vanishing of Ext modules.
\end{abstract}
\maketitle
\section{Introduction}

Let $R$ be a commutative Noetherian ring.
Goto \cite{G} proved the following theorem.

\begin{thm}[Goto]\label{goto}
The following are equivalent:
\begin{enumerate}
\item[{\rm (1)}]
$R$ is Gorenstein;
\item[{\rm (2)}]
For every finite $R$-module $M$, there exists an integer $n$ such that $\Ext _R^i(M,R)=0$ for all $i>n$.
\end{enumerate}
\end{thm}

We should note that this theorem remains valid even in the case where the ring $R$ has infinite Krull dimension.

The purpose of this note is to give a characterization of finite modules locally of finite injective dimension.
Our theorem is the following.

\begin{thm}\label{main}
The following are equivalent for a finite $R$-module $N$:
\begin{enumerate}
\item[{\rm (1)}]
$\id _{R_\p}\,N_\p <\infty$ for every $\p\in\Spec\,R$;
\item[{\rm (2)}]
For every finite $R$-module $M$, there exists an integer $n$ such that $\Ext _R^i(M,N)=0$ for all $i>n$.
\end{enumerate}
\end{thm}

This theorem is a generalization of Goto's.
In fact, applying our theorem to $N=R$, we immediately obtain Goto's theorem.

\section{Proof of the thoerem}

We denote by $\CM (R)$ the {\it Cohen-Macaulay locus} of $R$, that is, the set of prime ideals $\p$ of $R$ such that the local ring $R_\p$ is Cohen-Macaulay.
The following lemma can be shown in a similar way to the proof of \cite[Theorem 24.5]{M}.

\begin{lem}\label{nccm}
Let $\p$ be a prime ideal of $R$ such that both $R_\p$ and $R/\p$ are Cohen-Macaulay rings.
Then there exists an element $f\in R-\p$ such that $D(f)\cap V(\p )\subseteq \CM (R)$.
\end{lem}

It is known that the Cohen-Macaulay locus of a homomorphic image of a Cohen-Macaulay ring is an open subset; see \cite[Exercises 24.2]{M}.
To prove our theorem, we need to generalize this fact.
For an ideal $I$ of $R$, let $\CM _R(R/I)$ denote the set of prime ideals $\p\in V(I)$ such that the local ring $(R/I)_\p$ is Cohen-Macaulay.

\begin{lem}\label{relopen}
Let $I$ and $J$ be ideals of $R$.
Suppose that $V(J)$ is contained in $\CM (R)$.
Then the following hold:
\begin{enumerate}
\item[{\rm (1)}]
For any prime ideal $\p\in\CM _R(R/I)\cap V(J)$, there exists an element $f\in R-\p$ such that $D(f)\cap V(\p)\subseteq\CM _R(R/I)$.
\item[{\rm (2)}]
There exists an ideal $K$ of $R$ such that
$$
\CM _R(R/I)\cap V(J)=D(K)\cap V(I)\cap V(J).
$$
In other words, $\CM _R(R/I)\cap V(J)$ is an open subset of $V(I)\cap V(J)$ in the relative topology induced by the Zariski topology of $\Spec\,R$.
\end{enumerate}
\end{lem}

\begin{pf}
(1) Let $\p\in\CM _R(R/I)\cap V(J)$.
Then by the assumption that $V(J)$ is contained in $\CM (R)$, the ring $R_\p$ is a Cohen-Macaulay local ring.
Making a similar argument to the proof of \cite[Theorem 24.5]{M}, we can assume without loss of generality that there is an $R$-regular sequence $\xx =x_1,x_2,\dots, x_n$ in $\p$ with $\p ^r\subseteq\xx R$ for some $r>0$ and that $\overline\p ^i/\overline\p ^{i+1}$ is a free $\overline R/\overline\p$-module for all $i>0$, where $\overline R = R/\xx R$ and $\overline\p = \p /\xx R$.

We have only to prove that the residue ring $R/\p$ is Cohen-Macaulay.
In fact, if $R/\p$ is Cohen-Macaulay, then so are $(R/I)/(\p /I)$ and $(R/I)_{\p /I}$ since $\p$ is in $\CM _R(R/I)$.
Hence Lemma \ref{nccm} implies that there is an element $f\in R-\p$ such that $D(\overline f)\cap V(\p /I)$ is contained in $\CM (R/I)$, where $\overline f$ denotes the residue class of $f$ in $R/I$.
We easily see that $D(f)\cap V(\p)$ is contained in $\CM _R(R/I)$.

Let us show that $R/\p$ is a Cohen-Macaulay ring.
It is easy to see from \cite[Exercises 24.1]{M} that $R/\p =\overline R/\overline\p$ is Cohen-Macaulay if and only if so is $\overline R$.
Take a prime ideal $\q\in V(\xx)=V(\p)$.
Then we have $\q\supseteq\p\supseteq J$, hence $\q\in V(J)\subseteq\CM (R)$.
Therefore $R_\q$ is a Cohen-Macaulay local ring, and so is $\overline R_\q$ because $\xx$ is an $R_\q$-regular sequence.
This shows that $\overline R$ is a Cohen-Macaulay ring.
Thus we conclude that the residue ring $R/\p$ is Cohen-Macaulay, as desired.

(2) Set $U=\big\{\,\p /I+J\,\big |\,\p\in\CM _R(R/I)\cap V(J)\,\big\}$.
This is a subset of $\Spec\,R/I+J$.
Note that this subset is stable under generalization.
Let $P\in U$.
Then there is a prime ideal $\p\in\CM _R(R/I)\cap V(J)$ such that $P=\p /I+J$.
By the assertion (1) of the lemma, the set $D(f)\cap V(\p)$ is contained in $\CM _R(R/I)$ for some $f\in R-\p$.
Denote by $\overline f$ the residue class of $f$ in $R/I$.
It is easy to see that $P\in D(\overline f)\cap V(P)\subseteq U$.
Thus $U$ contains a nonempty open subset of $V(P)$.
By virtue of topological Nagata criterion \cite[Theorem 24.2]{M}, $U$ is an open subset of $\Spec\, R/I+J$; we have $U=D(K/I+J)$ for some ideal $K$ of $R$ containing $I+J$.
Then it is easily checked that $\CM _R(R/I)\cap V(J) = D(K)\cap V(I)\cap V(J)$.
\qed
\end{pf}

Now, we can prove our theorem.

\begin{tpf}
(2) $\Rightarrow$ (1): Let $\p$ be a prime ideal of $R$.
Then there is an integer $n$ such that $\Ext _R^i(R/\p , N)=0$ for all $i>n$.
Hence we have $\Ext _{R_\p}^i(\kappa (\p), N_\p)=0$ for all $i>n$.
Therefore by \cite[Theorem 3.1.14]{BH} we obtain $\id _{R_\p}\,N_\p\le n <\infty$.

(1) $\Rightarrow$ (2): First of all, note that (2) is equivalent to the statement that for each ideal $I$ of $R$ there is an integer $n$ such that $\Ext _R^i(R/I,N)=0$ for all $i>n$.
(This can easily be proved by induction on the number of generators of the $R$-module $M$.)
Suppose that there exists an ideal $I$ of $R$ such that for any integer $n$ there is an integer $i>n$ such that $\Ext _R^i(R/I,N)\ne 0$.
We want to derive a contradiction.
Since $R$ is Noetherian, one can choose $I$ to be a maximal one among such ideals.
Making a similar argument to the proof of Theorem \ref{goto}, we see that the ideal $I$ is prime and that for any element $f\in R-I$, the map
$$
\Ext _R^i(R/I,N)\overset{f}{\to}\Ext _R^i(R/I,N)
$$
is an isomorphism for $i\gg 0$.

\begin{claim}\label{supp}
One has $I\in\Supp _R\, N\subseteq\CM (R)$.
\end{claim}

\begin{cpf}
Our assumption (1) implies that for any $\p\in\Supp _R\,N$, the nonzero finite $R_\p$-module $N_\p$ has finite injective dimension.
Hence $R_\p$ is a Cohen-Macaulay local ring; see \cite[Corollary 9.6.2 and Remarks 9.6.4]{BH}.
Thus $\Supp _R\,N$ is contained in $\CM (R)$.
On the other hand, assume that $I$ is not in $\Supp _R\,N$.
Then there exists an element $f\in\Ann _R\,N-I$, and the map $\Ext _R^i(R/I,N)\overset{f}{\to}\Ext _R^i(R/I,N)$ is an isomorphism for $i\gg 0$.
Since $fN=0$, this map is the zero map, and we get $\Ext _R^i(R/I, N)=0$ for $i\gg 0$.
It follows from this contradiction that $I$ belongs to $\Supp _R\,N$.
\qed
\end{cpf}

Noting that $\Supp _R\,N=V(\Ann _R\,N)$, we see from Claim \ref{supp} and Lemma \ref{relopen}(2) that there is an ideal $K$ of $R$ such that $\CM _R(R/I)\cap\Supp _R\,N=D(K)\cap V(I)\cap\Supp _R\,N$.
The localization $(R/I)_I=\kappa (I)$ is a field, hence a Cohen-Macaulay ring.
It is seen from Claim \ref{supp} again that $I\in\CM _R(R/I)\cap\Supp _R\,N\subseteq D(K)$.
Thus there is an element $f\in K-I$.

\begin{claim}\label{is}
For any prime ideal $\p\in D(f)$ and any integer $i>\height\, I$, one has $\Ext _{R_\p} ^i(R_\p /IR_\p,N_\p)=0$.
\end{claim}

\begin{cpf}
We may assume that $\p$ belongs to both $V(I)$ and $\Supp _R\,N$ because otherwise the module $\Ext _{R_\p} ^i(R_\p /IR_\p,N_\p)$ automatically vanishes.
Hence Claim \ref{supp} implies that $\p$ belongs to $\CM (R)$, namely the local ring $R_\p$ is Cohen-Macaulay.
Added to it, since $D(f)$ is contained in $D(K)$, we have $\p\in D(K)\cap V(I)\cap \Supp _R\,N\subseteq \CM _R(R/I)$, and therefore $R_\p /IR_\p$ is Cohen-Macaulay.
Thus we get the following equalities:
$$
\depth\,R_\p - \depth\,R_\p /IR_\p = \dim R_\p - \dim R_\p /IR_\p = \height\, IR_\p = \height\, I.
$$
Since $N_\p$ is a finite $R_\p$-module of finite injective dimension by assumption, it follows from the result of Ischebeck \cite[Exercises 3.1.24]{BH} that $\Ext _{R_\p} ^i(R_\p /IR_\p,N_\p)=0$ for every $i>\height\,I$.
\qed
\end{cpf}

Claim \ref{is} means that $(\Ext _{R_f}^i(R_f/IR_f, N_f))_\P =0$ for every $\P\in\Spec\,R_f$ and every $i>\height\, I$.
Therefore $\Ext _{R_f}^i(R_f/IR_f, N_f)=0$ for $i>\height\, I$.
The $R$-module $\Ext _R ^i(R/I,N)$ is isomorphic to $\Ext _{R_f}^i(R_f/IR_f, N_f)$ for $i\gg 0$, and thus $\Ext _R^i(R/I,N)=0$ for $i\gg 0$.
This contradiction completes the proof of our theorem.
\qed
\end{tpf}



\end{document}